\documentclass{amsart}

\usepackage[T1]{fontenc}
\usepackage{mathtools}
\usepackage{tensor} 
\usepackage{amssymb} 
\usepackage{enumitem} 
\usepackage{accents}
\usepackage{mathrsfs} 
\newcommand{\ubar}[1]{\underaccent{\bar}{#1}} 

\newtheorem{theorem}[subsection]{Theorem}
\newtheorem{lemma}[subsection]{Lemma}
\newtheorem{proposition}[subsection]{Proposition}
\newtheorem{corollary}[subsection]{Corollary}
\newtheorem{definition}[subsection]{Definition}
\newtheorem{remark}[subsection]{Remark}
\allowdisplaybreaks

\title{CONFORMAL CIRCLES AND LOCAL DIFFEOMORPHISMS}
\author{Tzu-Mo Kuo}
\address{Department of Mathematics, University of California, Santa Cruz, CA 95064, USA}
\email{tkuo6@ucsc.edu}

\begin{document}
\maketitle

\begin{abstract}
We study unparametrized conformal circles, or called conformal geodesics, study diffeomorphisms mapping conformal circles to conformal circles in pseudo-Riemannian conformal manifolds. We show that such local diffeomorphisms are conformal local diffeomorphisms. Our result extends the result of Yano and Tomonaga. We also present a holographic interpretation for our result on Poincaré-Einstein manifolds. The proofs take suitable variations of conformal circles.
\end{abstract}

\section{Introduction}
\label{sec:intro}
Riemannian geodesics, as fundamental geometric objects, are often considered when studying Riemannian structures. One classic problem concerning geodesics and Riemannian structures is the following: If there is a diffeomorphism that maps geodesics to geodesics, is it an isometry? The answer is negative due to affine transformations on Euclidean spaces. In general, one may need to further assume irreducible Riemannian manifolds for the diffeomorphism to be an isometry \cite{irredriemannaffinetransf,foundationtodiffgeom}. The parallel problems for CR manifolds \cite{Chainpreserving} and conformal manifolds \cite{Caratheodory,preserveconformalcircleinf} are affirmative in some sense.

In the context of a pseudo-Riemannian conformal manifold $(M^{n},[g])$ with $n\geq 2$, a distinguished family of curves known as conformal circles or conformal geodesics emerges. These curves satisfy a third-order differential equation for non-null circles \cite{BaileyandEastwood}. The derivation of these conformal circles is based on various perspectives of conformal manifolds, including Cartan geometry \cite{FriedrichandSchmidt}, the standard tractor bundle \cite{tractorbundle}, and the Poincaré-Einstein manifold \cite{JoelfineandHerfray}. More detailed discussions about each viewpoint on conformal manifolds can be found in references such as \cite{tractorbundlelecture, FeffermanandGraham2,geometrytextbook3,tractorambient, parabolictextbook}. Taking the unit sphere $(S^n,[g_{S^n}])$ as a Riemannian conformal model, the conformal circles in the model are either straight lines or planar circles when viewed through stereographic projection onto Euclidean space \cite{Tod_exampleofconformalcircle}. In the context of the Poincaré ball $B^{n+1}$, each circle on the boundary $S^n=\partial B^{n+1}$ can be orthogonally extended to form a totally geodesic surface within $B^{n+1}$.  Based on the case of the conformally flat model, each conformal circle in a Riemannian conformal manifold $(M,[g])$ can be formally extended to an asymptotically totally geodesic surface in the Poincaré-Einstein space $(M_+,g_+)$ by holographic construction \cite{JoelfineandHerfray} where the term holography comes from physics, e.g., \cite{Witten1998AntideSS}. We will use the term "holography" to mean the geometry from Poincaré-Einstein space $(M_+,g_+)$. 

In this article, we consider the parallel classic problem for conformal circles with a local diffeomorphism $f$ between pseudo-Riemannian conformal manifolds $(M^{n},[g])$, $(N^n,[h])$ with the same signature $(p,q)$ for $n\geq 2$ and the problem in holographic settings $F\colon M_+\to N_+$. The problem for conformal circles in a Riemannian conformal manifold can be traced back to Carathéodory \cite{Caratheodory}. He considered a bijection on $\mathbb R^2$, which doesn't need to be continuous, that maps straight lines (resp. circles) to straight lines (resp. circles), and he showed the bijection is a conformal transformation. Later on, K. Yano and Y. Tomonaga \cite{Yano, preserveconformalcircleinf} showed that an infinitesimal transformation is a conformal killing vector field if and only if it carries unit-speed conformal circles to unit-speed conformal circles. The equivalence problem is also discussed in Cartan geometry settings in terms of distinguished curves \cite{Cartangeomdistcurvequivproblem}. In \textsection \ref{sec:confcirc}, we review conformal circles and prove Theorem \ref{thm:conformalinftycircle} for the equivalence problem of conformal circles. Our proof is motivated by the work from Yano and Tomonaga. In \textsection\ref{sec:holointerp}, we review Poincaré-Einstein manifolds and extend the results from Fine and Herfray \cite{JoelfineandHerfray} to indefinite signature. In \textsection\ref{sec:pfofthm:holocirc}, we prove Theorem \ref{thm:holocircle2} in holographic settings.\\

The author would like to thank his advisor, Jie Qing, for all discussions and illuminating advice.

\section{Conformal Circles on Conformal Manifolds}\label{sec:confcirc}
In this section, we review definitions of conformal circles on a pseudo-Riemannian conformal manifold $(M^n,[g])$. As mentioned in the introduction, we prove Theorem \ref{thm:conformalinftycircle} which extends the result of Yano and Tomonaga \cite{Yano}. 

Let $\hat g=e^{2\sigma} g$ for $g\in [g]$ and $\sigma\in C^{\infty}(M)$. The Levi-Civita connections $\widehat\nabla$ and $\nabla$ of respective $\hat g$ and $g$ on sections of $TM$ are related by \cite{LubbeandTod}
\begin{align}\label{eq:transfofconn}
    \widehat\nabla_i v^j=\nabla_iv^j+S\indices{_i_k^j^l}v^k\nabla_l\sigma
\end{align}
where the Latin indices $i,$ $j,$ $k,$ $l$ follow the convention in \cite{penroseindex}. The $S\indices{_i_k^j^l}$ is a conformally invariant tensor defined by 
\begin{align}
    S\indices{_i_k^j^l}=\delta^j_i\delta^l_k+\delta^l_i\delta^j_k-g_{ik}g^{jl}
\end{align}
where $\delta^j_i=g^{jk}g_{ki}$. We denote the Schouten tensor of $g$ and $\hat g$ by $P_{ij}$ and $\hat P_{ij}$ respectively for $n\geq 3.$ Their relations are
\begin{align}
 \label{eq:conftranofP}
 \hat P_{ij}=P_{ij}-\nabla_i\nabla_j\sigma+\nabla_i\sigma\nabla_j\sigma-\frac{1}{2}g_{ij}g^{kl}\nabla_k\sigma\nabla_l\sigma. 
\end{align}
When $n=2$, we assume that (\ref{eq:conftranofP}) is a part of the structure on $(M,[g]).$ \\

Recall that Riemannian geodesics are the projections of integral curves of constant horizontal vector fields on linear frame bundle (\cite{foundationtodiffgeom} Chapter 3 Proposition 6.3); conformal circles can be defined in a similar way from Cartan geometry \cite{FriedrichandSchmidt}.
    \begin{definition}\label{def:conformalcircle}
    Let $g\in [g]$ and $\nabla$ be the corresponding Levi-Civita connection. A conformal circle with respect to $g$ is a curve $\gamma\colon I\to M$ and a 1-form $b$ along $\gamma$ such that they satisfy
    \begin{align}
    \label{eq:confcirceq1}   \nabla_{\dot\gamma}\dot\gamma^i&=-S\indices{_j_k^i^l}b_l\dot\gamma^j\dot\gamma^k,\\
    \label{eq:confcirceq2}    \nabla_{\dot\gamma}b_i&=(\textstyle\frac{1}{2}b_jb_lS\indices{_k_i^j^l}+P_{ki})\dot\gamma^k
    \end{align}
    where the overdot of $\gamma$ denotes the derivative with respect to $t\in I$. The equations are conformally invariant for conformal change $\hat g=\Omega^2g$ if $\dot{\hat \gamma}=\dot\gamma$ and $\hat b_i= b_i-\Omega^{-1}\nabla_i\Omega$. 
    \end{definition}
Denote $g_{ij}v^iw^j$ by $\langle v,w\rangle$ for $g\in [g]$. We call a curve $\gamma\colon I\to M$ \textit{null} if $\langle \dot\gamma,\dot\gamma\rangle =0$ on $I$ and it is called \textit{non-null} if it's not null. Direct computation from (\ref{eq:confcirceq1}) shows
\begin{equation}
    \label{eq:nullityofconfcirc}
    \nabla_{\dot\gamma}\langle\dot\gamma,\dot\gamma\rangle=-2\langle\dot\gamma,\dot\gamma\rangle b_i\dot\gamma^i.
\end{equation}
Therefore, if a conformal circle has a null velocity at some point, then it's a null conformal circle. For completeness, we recall null pseudo-Riemannian geodesics are necessary and sufficient to be null conformal circles with some reparametrization \cite{Lubbe,Tod_exampleofconformalcircle}. 

Assume a conformal circle $\gamma$ is non-null. By solving $b_i$ in (\ref{eq:confcirceq1})
\begin{equation}\label{eq:bform}
    b_i=\frac{1}{\langle\dot\gamma,\dot\gamma\rangle}\left(\nabla_{\dot\gamma}\dot\gamma_i-2\frac{\langle \dot\gamma,\nabla_{\dot\gamma}\dot\gamma\rangle}{\langle\dot\gamma,\dot\gamma\rangle}\dot\gamma_i\right),
\end{equation} 
one can have a third-order differential equation from (\ref{eq:confcirceq2}). The third-order differential equation is equivalent to the system of the equations (\ref{eq:confcirceq1}) and (\ref{eq:confcirceq2}) if one defines the one form $b_i$ back.
\begin{definition}\label{def:nonnullcircle}
Given $g\in[g]$ with corresponding Levi-Civita connection $\nabla$. A parametrized non-null conformal circle $\gamma$ is defined to satisfy
\begin{align}\nabla_{\dot\gamma}\nabla_{\dot\gamma}\dot\gamma^i=3\frac{\langle\dot\gamma, \nabla_{\dot\gamma}\dot\gamma\rangle}{\langle\dot\gamma,\dot\gamma\rangle}\nabla_{\dot\gamma}\dot\gamma^i-\frac{3 \langle\nabla_{\dot\gamma}\dot\gamma,\nabla_{\dot\gamma}\dot\gamma\rangle}{2\langle \dot\gamma,\dot\gamma\rangle}\dot\gamma^i+\langle \dot\gamma,\dot\gamma\rangle\dot\gamma^jP\indices{_j^i}-2P_{jk}\dot\gamma^j\dot\gamma^k\dot\gamma^i \label{eq:nonnullconformalcircle}
\end{align}
with initial conditions $\gamma(0),\;\dot\gamma(0),\;\nabla_{\dot\gamma}\dot\gamma(0)$. The equation is invariant under conformal change $\hat g=\Omega^2 g$.
\end{definition}
Since the induced metric on a non-null curve $\gamma$ is nondegenerate, we can consider the orthogonal decomposition of the pullback bundle of $TM$ by $\gamma$ \cite{semiriem}. We call $\gamma$ satisfies the \textit{tangential (resp. normal) part} of (\ref{eq:nonnullconformalcircle}) if it is a solution of the equation that is the orthogonal projection of (\ref{eq:nonnullconformalcircle}) to the tangent (resp. normal) bundle of $\gamma$. It is known \cite{BaileyandEastwood} that any regular curve can be reparametrized to satisfy the tangential part of (\ref{eq:nonnullconformalcircle}). The normal part of (\ref{eq:nonnullconformalcircle}) is invariant under reparametrization of $\gamma$ and it is only satisfied by non-null conformal circles. Since (\ref{eq:nonnullconformalcircle}) is derived from (\ref{eq:confcirceq2}), it's convenient to introduce a vector field along an arbitrary curve $\gamma$ 
\begin{equation}\label{eq:vectfromconfcirceq2}
    E^i(\gamma,v,g)=\nabla_{\dot\gamma}v^i-(\textstyle\frac{1}{2}v^jv^lS\indices{_k^i_j_l}+P\indices{_k^i})\dot\gamma^k
\end{equation}
where  $v$ is a vector field along $\gamma$. If $v$ satisfies the right-hand side of (\ref{eq:bform}) by descending index, then we denote the vector field $E^i(\gamma,v,g)$ by $E^i(\gamma,g).$

For completeness, we recall another third-order differential equation for a non-null conformal circle $\gamma$. If $\gamma\colon I\to M$ is reparametrized so that it is of unit tangent velocity with respect to $g\in[g]$, then it satisfies \cite{Tod_exampleofconformalcircle}
\begin{align*}
\nabla_{\dot\gamma}\nabla_{\dot\gamma}\dot\gamma^i&=-\left(\langle\nabla_{\dot\gamma}\dot\gamma,\nabla_{\dot\gamma}\dot\gamma\rangle+P_{jk}\dot\gamma^j\dot\gamma^k\right)\dot\gamma^i+P^i_j\dot\gamma^j \quad\text{if }\langle\dot\gamma,\dot\gamma\rangle =1;\\
\nabla_{\dot\gamma}\nabla_{\dot\gamma}\dot\gamma^i
&=\left(\langle\nabla_{\dot\gamma}\dot\gamma,\nabla_{\dot\gamma}\dot\gamma\rangle-P_{ik}\dot\gamma^i\dot\gamma^k\right)\dot\gamma^j+P^j_i\dot\gamma^i \qquad\text{if }\langle\dot\gamma,\dot\gamma\rangle =-1;
\end{align*}
The first line is the equation introduced by Yano \cite{Yano}.
\\

Given a local diffeomorphism $f\colon M\to N$ between pseudo-Riemannian conformal manifolds $(M^n,[g])$ and $(N^n,[h])$. Assume both of the conformal classes have same signature $(p,q)$. If $f$ is a conformal local diffeomorphism, it's direct to see $f$ maps unparametrized conformal circles to unparametrized conformal circles. The converse direction is also true if the map $f$ preserves some nullity condition.
\begin{theorem} \label{thm:conformalinftycircle}
Let $(M^n,[g])$ and $(N^n,[h])$ be pseudo-Riemannian conformal manifolds with same signature $(p,q)$. Assume a local diffeomorphism $f\colon M\to N$ satisfying 
\begin{enumerate}[label=(\roman*)]
    \item $df(v)$ is non-null (resp. null) $\;\forall$non-null (resp. null) $v\in TM$ if $p\neq q$\label{statement-1 in thm of conformal circle};
    \item sgn$\langle v,v\rangle_g$=\;sgn$\langle df(v),df(v)\rangle_h$ $\forall v\in TM$ if $p=q$.\label{statement-2 in thm of conformal circle}
\end{enumerate}
If $f$ maps unparametrized non-null conformal circles to unparametrized non-null conformal circles, then $f$ is a conformal local diffeomorphism.
\end{theorem}
\begin{proof}
Given $y\in M$. Let $g\in[g]$ and $h\in[h]$. Choose a normal coordinate of $g$ centered at $y$, $\{x^i|\,1\leq i\leq n\}$. Since $f$ is a local diffeomorphism, we can identify the coordinate system near $f(y)$ as $\{x^i\}_{i=1}^n$. Let $\gamma(t)$ be a parametrized non-null conformal circle satisfying (\ref{eq:nonnullconformalcircle}) with initial conditions $\dot\gamma^k_0$, $\ddot\gamma^k_0$ at $y=\gamma(0)$ where $\ddot\gamma^k_0$ is the coordinate of $\nabla_{\dot\gamma}\dot\gamma(0)$ with respect to $g$. Since $f\circ\gamma$ is an unparametrized non-null conformal circle, it satisfies the normal part of (\ref{eq:nonnullconformalcircle}) with the given parameter $t$, which is the following in the coordinate we chose.
\begin{align}\label{eq:normalptofnonnullconfcirceq}
E^k(\gamma, h)-\frac{\langle E(\gamma, h),\dot\gamma\rangle_h}{\langle \dot\gamma,\dot\gamma\rangle_h}\dot\gamma^k=0
\end{align}
where $E(\gamma, h)$ is a vector field along $\gamma$ defined from (\ref{eq:vectfromconfcirceq2}). In the following, we are considering $t=0$ for (\ref{eq:normalptofnonnullconfcirceq}). Observe that (\ref{eq:normalptofnonnullconfcirceq}) is a degree-two polynomial of $\ddot\gamma^k_0$ with coefficients depending on the derivatives of $g$ and $h$. It is because $\gamma(t)$ satisfies the third order differential equation (\ref{eq:nonnullconformalcircle}) with respect to $g$. Now let $\dot\gamma^k_0=V^k$ and $\ddot\gamma^k_0=\epsilon A^k$ where $V,A\in \mathbb R^n$ are fixed. The variable $\epsilon\in\mathbb R$ is an arbitrary number in an open interval containing $1$. Because the non-null conformal circle equation is an autonomous ODE, (\ref{eq:normalptofnonnullconfcirceq}) depends smoothly on $\epsilon$. Based on the arguments we just made, we know (\ref{eq:normalptofnonnullconfcirceq}) at $t=0$ is a polynomial of $\epsilon$ with degree two. Therefore, the coefficient of $\epsilon^2$ vanishes. After direct computation, the coefficient gives
\begin{align}
\frac{\langle V, A\rangle_g}{\langle V, V\rangle_g} \left(A^k-\frac{\langle V,A\rangle_{h}}{\langle V,V\rangle_{h}}V^k\right)
=
\frac{\langle V, A\rangle_{h}}{\langle V,  V\rangle_{h}}\,A^k
-
\frac{\langle  V, A\rangle^2_{h}}{\langle  V, V\rangle^2_{h}}\,V^k.
\end{align}
If $\langle V,A\rangle_g =0$, we get
\begin{align}\label{eq:coefficientofxi^2}
    0=\frac{\langle V, A\rangle_{h}}{\langle V, V\rangle_{h}}\,A^k
    -
    \frac{\langle  V, A\rangle^2_{h}}{\langle V, V\rangle^2_{h}}\,V^k.
\end{align}
Assume the normal coordinate we chose is $g_{ii}(y)>0$ for $1\leq i\leq p$ and $g_{ii}(y) <0$ for $p+1\leq i\leq p+q$. If $V^k=\delta^{ik}$ and $A^k=\delta^{jk}$ for $i\neq j$, then $h_{ij}(y)=0$ from (\ref{eq:coefficientofxi^2}). If $g_{ii}=g_{jj}$ for $i\neq j$, we then let $V^k=\delta^{ik}+\delta^{jk}$ and $A^k=\delta^{ik}-\delta^{jk}$; so, we have $h_{ii}=h_{jj}$ from (\ref{eq:coefficientofxi^2}). Therefore, the pullback metric $(f^*h)_{ij}$ is of the form 
$$\begin{pmatrix}
B\mathbb{I}_p &{}\\
{ }& -C\mathbb{I}_q
\end{pmatrix} $$
for some $B,\,C\neq 0$. If $p\neq q$, then $B$ and $C$ are positive because $h$ is of the signature $(p,q)$. If $p=q$, then $B$ and $C$ can be both positive or both negative. Recalling that the $f$ preserves the nullity of null vectors, we know $(\partial_i+\partial_j)$ is null at $y$ with respect to $f^*h$ when $g_{ii}\neq g_{jj}$ which implies $B=C$. The sign of $B$ and $C$ is positive for $p=q$ since the sign of $h_{ii}$ is the same as the sign of $g_{ii}$.
\end{proof}

\begin{remark}
    Note that if $(M,[g])$ is Riemannian, that is $p=n$ and $q=0$, then any local diffeomorphism $f\colon M\to N$ automatically satisfies \ref{statement-1 in thm of conformal circle} and \ref{statement-2 in thm of conformal circle} in Theorem \ref{thm:conformalinftycircle}.
\end{remark}
\begin{remark}
The particular case for Theorem \ref{thm:conformalinftycircle} has been studied in the literature if one assumes the $f$ to be a bijection (no need to be continuous) and $M=N=\mathbb R^n$ with the standard Riemannian conformal structure \cite{Carathéodory, J.Jeffers}. Note that the proof in \cite{Carathéodory,J.Jeffers} needs the global property of conformal circles; that is, $f$ maps straight lines (resp. circles) to straight lines (resp. circles). However, by assuming additional regularity of $f$ in this paper, we only need the local condition of conformal circles, namely the conformal circle equation, to establish Theorem \ref{thm:conformalinftycircle}.   
\end{remark}
\begin{remark}
    One can give a different proof from (\cite{Spivak2}, Chapter 6, Addendum 1) for the parallel problem of parametrized Riemannian geodesics by following the proof idea of Theorem \ref{thm:conformalinftycircle}.
\end{remark}
\begin{remark}\label{rmk:pfofholocircle}
    Though the initial condition $\ddot \gamma^k_0=\epsilon A^k$ in the proof gives an $\epsilon$-family of conformal circles which induce a Jacobi field, we do not need the Jacobi field equation introduced by \cite{Friedrichconformalgeod,LubbeandTod} to prove the theorem. 
\end{remark}

We follow similar arguments of Theorem \ref{thm:conformalinftycircle} to prove Theorem \ref{thm:holocircle2} in \textsection \ref{sec:pfofthm:holocirc}.\\

\section{Holography Interpretation for Conformal Circles}
\label{sec:holointerp}
In this section, we review Poincaré-Einstein manifolds and extend the results from Fine and Herfray \cite{JoelfineandHerfray} to indefinite signature.

\subsection*{Poincaré-Einstein manifold}
\hfill

Let $(M^n,[g])$ be a pseudo-Riemannian conformal manifold with signature $(p,q)$ and $n\geq 2$. There exists a pseudo-Riemannian manifold $(M^{n+1}_+,g_+)$ with boundary $\partial M_+=M$ \cite{FeffermanandGraham2}. The signature of $g_+$ is $(p+1,q)$. Let $r\in C^\infty (M_+)$ be a defining function for $M$, that is $r>0$ on the interior $M^0_+$, $r=0$ on $M$ and $dr\neq 0$ on $M$.  Then, $g_+$ and $r$ satisfy
\begin{enumerate}[label=(\roman*)]
    \item $r^2g_+$ can be smoothly extended to be a metric on $M_+$ so that $$\bar g|_M\triangleq r^2g_+|_M\in [g];$$
    \item $Ric(g_+)+ng_+=O(r)$.
\end{enumerate}
The pair $(M_+,g_+)$ is called a \textit{Poincaré-Einstein manifold} for $(M,[g])$ and the pair $(M,[g])$ is called the \textit{conformal infinity}  to $(M_+,g_+)$. 

Let $g\in[g]$. There exists a unique defining function $r$, called \textit{geodesic defining function}, such that $|dr|^2_{\bar g}=1$ near $M\subset M_+$. The function $r$ makes an identification between a neighborhood $U$ of $M$ in $ M_+$ and a neighborhood $\mathcal U$ of $M\times\{0\}$ in $ M\times [0,\infty)$. By the identification, $g_+$ is in \textit{normal form} relative to $g$; that is, $g_+=\frac{dr^2+g_r}{r^2}$ on $\mathcal U$ where $g_r=g-r^2P+O(r^3)$. The tensor $P$ is the Schouten tensor of $g$ for $n\geq 3$. When $n=2$, $P$ is a symmetric two-tensor on $M$ satisfying $P\indices{_i^i}=\frac{1}{2}R$ and $P\indices{_i_j_,^j}=\frac{1}{2}R_{,i}$ where $R$ is the scalar curvature of $g$. By the pulling back of an even diffeomorphism between neighborhoods of $M\times\{0\}\subset M\times [0,\infty)$ which restricts to the identity map on $M\times\{0\}$, the normal forms $g_+$ relative to conformal related metrics are identical modulo $O(r)$. 
\begin{remark}
    Note that the orders of $r$ above in the Ricci condition and in $g_r$ can be further refined to higher orders depending on dimension $n$. See details in \cite{FeffermanandGraham2}.
\end{remark}
\begin{remark}
    For $n=2$, the trace and the divergence conditions of $P$ are conformally invariant. ( \cite{FeffermanandGraham2}, arguments after Theorem 3.7)
\end{remark}

\subsection*{Surfaces in the Poincaré-Einstein Manifold}
\hfill

Let $(M^{n+1}_+,g_+)$ be a Poincaré-Einstein manifold in normal form relative to $g\in [g]$ where $r$ is its corresponding geodesic defining function. Let $\gamma\colon I\to  M$ be a non-null curve. The interval $I$ can be shrunk if necessary.  Choose a local coordinate $\{x^i|\, 1\leq i\leq n\}$ on an open set $\mathcal W$ in $M$ containing $\gamma(t)$ for $t\in I$. The coordinate of $\gamma(t)$ is denoted by $\gamma^i(t)$. If $\Sigma\subset \overline{M_+}$ is an embedded surface orthogonal to $M$ with $\Sigma\cap M=\gamma$, then one can have an \textit{asymptotic isothermal coordinate} of $\Sigma$ near $\gamma$; that is, there is a diffeomorphism $\sigma\colon (t,\lambda)\mapsto (x^i(t,\lambda),r(t,\lambda))$ from $I\times I$ to $\Sigma\subset \overline{M_+}$ such that
    \begin{equation}
    \label{eq:condiforasymptisothermalcoord}
        \begin{cases}
    \sigma(t,0)=\gamma^i(t)
    \\
    \sigma^*\bar g=\begin{pmatrix}(-1)^\varepsilon &0\\
    0 &1\end{pmatrix}c(t,\lambda)+O(\lambda^3)
    \end{cases},
    \end{equation}
where $c(t,0)\neq 0$ and $\varepsilon=0$ if $\langle\dot\gamma,\dot\gamma\rangle >0$ and $\varepsilon=1$ otherwise. In fact, to satisfy (\ref{eq:condiforasymptisothermalcoord}) for $\Sigma$ orthogonal to $M$, the expansions of $x^i$, $r$ and $\sigma^*g_+$ with respect to $\lambda$ are in the following forms.
    \begin{proposition}\cite{JoelfineandHerfray}\label{prop:expofrandxinisothermal}
    Let $g_+$ be a Poincaré-Einstein metric and $\Sigma\subset \overline{M_+}$ a surface as above. Then, the asymptotic isothermal coordinate in (\ref{eq:condiforasymptisothermalcoord}) satisfies
    \begin{align}
        \label{eq:expofrandx}
        &\small\begin{matrix*}[c]
         x^i(t,\lambda) &\equiv& \gamma^i(t)& +& 0  &+& \frac{|\dot\gamma|^2v^i}{2} \lambda^2 &+& \frac{u^i}{3}\lambda^3 
            \\
         r(t,\lambda) &\equiv&0   &+& |\dot\gamma|\lambda &+& 0 &+& \frac{(-1)^\varepsilon|\dot\gamma|}{6}\left[\kappa(\gamma,v,g)-\frac{3}{2}\langle\dot\gamma,\dot\gamma\rangle\langle v,v\rangle\right]\lambda^3
        \end{matrix*}
    \end{align}
    and
    \begin{equation}
    \label{eq:isothermalg+}\sigma^* g_+=\frac{1}{\lambda^2}\begin{pmatrix}(-1)^\varepsilon &0\\
        0 &1\end{pmatrix}
        \left(1+\frac{(-1)^\varepsilon\; 2}{3}\kappa(\gamma,v,g)\lambda^2\right)+O(\lambda),
    \end{equation}  
 where the expansions of $x^i$ and $r$ are modulo $O(\lambda^4)$, the $|\dot\gamma|$ denotes the square root of $|\langle\dot\gamma,\dot\gamma\rangle|$, the $v^i$ satisfies $\langle\dot\gamma,v\rangle=\left\langle\dot\gamma,\nabla_{\dot\gamma}\left(\frac{\dot\gamma}{\langle\dot\gamma,\dot\gamma\rangle}\right)\right\rangle$, the $u^i$ satisfies $\langle\dot\gamma,u\rangle=0$ and $\kappa(\gamma,v,g)=\langle E(\gamma,v,g),\dot\gamma\rangle$ which is the tangential part of (\ref{eq:vectfromconfcirceq2}). Note that $\nabla_{\dot\gamma}\left(\frac{\dot\gamma}{\langle\dot\gamma,\dot\gamma\rangle}\right)$ is equal to the right-hand side of the $b$-form in (\ref{eq:bform}) by lowering index.
\end{proposition}

 Since $\sigma^*\bar g$ is pseudo-Riemannian for small $\lambda$, the tangent bundle $TM_+$ has the orthogonal decomposition along $T\Sigma$ near $M$ \cite{semiriem}. Let $(e_{\alpha'})$ be a local orthonormal frame of the normal bundle of $\Sigma$ with respect to $\bar g$ near $\Sigma\cap M$. After direct computation, the Taylor expansions of $e_{\alpha'}$ with respect to $\lambda$ are \cite{JoelfineandHerfray}
\begin{align}
\label{eq:expofnormalvec}
e_{\alpha'}(t,\lambda)=&\phi_{\alpha'}(\lambda)-\frac{(-1)^\varepsilon\lambda^2}{2}\left[\langle \phi_{\alpha'}(\lambda),\partial_t v\rangle_g+\left(v^i\partial_ig-2P^g\right)(\phi_{\alpha'}(\lambda),\dot\gamma)\right]\dot\gamma
\\
\notag&-\langle \phi_{\alpha'}(\lambda),|\dot\gamma|v\rangle_g\;\lambda\;\partial_r+O(\lambda^3),    
\end{align}
where $\phi_{\alpha'}(\lambda)$ is a family of sections of the normal bundle of $\gamma $ in $M$.
\begin{proposition}\cite{JoelfineandHerfray}
    \label{prop:expof2ndfundform}
    Considering the projection of the second fundamental form of $\Sigma$ on $e_{\alpha'}$ with respect to $g_+$, then, its asymptotic expansion is
    \begin{equation}
    \label{eq:2ndformofsigma}
    \frac{1}{\lambda^2}
    \begin{pmatrix}
    (-1)^\varepsilon\left\langle\nabla_{\dot\gamma}\left(\frac{\dot\gamma}{\langle\dot\gamma,\dot\gamma\rangle}\right)- v,\phi_{\alpha'}(\lambda)\right\rangle-\frac{\langle u,\phi_{\alpha'}(\lambda)\rangle}{\langle\dot\gamma,\dot\gamma\rangle}\lambda  &\qquad*\\
    \langle E(\gamma,v,g),\phi_{\alpha'}(\lambda)\rangle\lambda  &\qquad(-1)^\varepsilon\frac{\langle u,\phi_{\alpha'}(\lambda)\rangle}{\langle\dot\gamma,\dot\gamma\rangle}\lambda
    \end{pmatrix}+O(1).
    \end{equation}
    Due to Proposition \ref{prop:expofrandxinisothermal}, the asymptotic minimal condition $H=O(r^2)$ of $\Sigma$ is equivalent to $v=\nabla_{\dot\gamma}\left(\frac{\dot\gamma}{\langle\dot\gamma,\dot\gamma\rangle}\right)$ which is exactly the same as the $b$-form in ($\ref{eq:bform}$). The asymptotic totally geodesic condition $K=O(r^2)$ is equivalently satisfied when $v=\nabla_{\dot\gamma}\left(\frac{\dot\gamma}{\langle\dot\gamma,\dot\gamma\rangle}\right)$, $u=0$ and $\gamma$ being an unparametrized conformal circle. 
\end{proposition}
 
\begin{definition}\label{def:propersurf}
    Let $\Sigma\subset \overline{M_+}$ be an embedded surface orthogonal to $M$ so that $\Sigma\cap M$ is a non-null curve. It is called a proper surface if it is asymptotic totally geodesic $K=O(r^2)$ where $K$ is the second fundamental form of $\Sigma$ with respect to $g_+$.
\end{definition}

\section{Proof of Theorem}
\label{sec:pfofthm:holocirc}
In this section, we consider a local diffeomorphism $F\colon M_+\to N_+$ which smoothly extends a local diffeomorphism $f\colon M\to N$. We introduce the definition of asymptotic local isometry and cosider its local conditions. We also introduce an adapted coordinate $\Sigma$ for the proof of Theorem \ref{thm:holocircle2}.

Let $(N^n,[h])$ be a pseudo-Riemannian conformal manifold with same signature $(p,q)$ as $(M,[g])$ and with a Poincaré-Einstein space $(N_+,h_+)$. We keep $\{x^i\}_{i=1}^n$ as a coordinate system on an open set $\mathcal W$ in $M$.
\begin{definition}\label{def:asymlocalisometry}
A local diffeomorphism $F\colon M_+\to N_+$ is called an asymptotic local isometry if 
\begin{align}\label{eq:asympisometry}
    F^*h_+-g_+=O(r).
\end{align}
\end{definition}

It's useful to realize Definition \ref{def:asymlocalisometry} in terms of local coordinates. Let $r$ and $s$ be geodesic defining functions for $g\in[g]$ and $h\in[h]$ respectively. Identifying some neighborhoods of $M\subset \overline{M_+}$ and $N\subset \overline{N_+}$ to neighborhoods $\mathcal U$ of $M\times\{0\}\subset M\times [0,\infty)$ and $\mathcal V$ of $N\times\{0\}\subset N\times [0,\infty)$ respectively, a local diffeomorphism $F\colon M_+\to N_+$ can be identified near $M\subset\overline{M_+}$ and $N\subset\overline{N_+}$ as a local diffeomorphism from $\mathcal U$ to $\mathcal V$ 
\begin{equation}
\label{eq:coorddescripofF}    \left(x,r\right)\mapsto\left(\mathscr F(x,r),F^s(x,r)\right)
\end{equation}
where $\mathscr F(x,0)=f(x)$, $s\circ F=F^s$ on $\mathcal U$ and $F^s(x,0)=0$. Then, (\ref{eq:asympisometry}) is equivalent to 
\begin{align}
    \label{eq:equivrelationofasympiso}
    F^*\bar h-\left(F^s/r\right)^2\bar g=O(r^3),
\end{align}
where $\bar g=dr^2+g_r$ and $\bar h=ds^2+h_s$. Since $f$ is a local diffeomorphism, we denote the coordinate of $\mathscr F(x,r)$ by $F^i(x,r)$ with $F^i(x,0)=x^i$. In terms of the coordinates $(x^i,r)$ on $M_+$ and $(x^i,s)$ on $N_+$, (\ref{eq:equivrelationofasympiso}) is given by
\begin{align}
O(r^3)&=F^s_{,i}F^s_{,j}+F^k_{,i}F^l_{,j}(h_s\circ F)_{kl}-\left(F^s/r\right)^2 (g_r)_{ij},
 \label{eq:1asympisomcondi}
 \\
 O(r^3)&=F^s_{,r}F^s_{,r}+F^k_{,r}F^l_{,r}\,(h_s\circ F)_{kl} - \left(F^s/r\right)^2,
 \label{eq:2asympisomcondi}
 \\
 O(r^3)&=F^s_{,r}F^s_{,i}+F^k_{,r}F^l_{,i}(h_s\circ F)_{kl}
 \label{eq:3asympisomcondi}
\end{align}
where we have used commas to express partial derivatives with respect to the coordinates $(x^i,r)$ on $M_+$. At $r=0$, (\ref{eq:1asympisomcondi}) and (\ref{eq:3asympisomcondi}) give
\begin{equation}
\label{eq:1stderivativeconditionforF}
\begin{split}
&f^*h=e^{2\sigma}g,\quad F^{s}_{,r}=e^\sigma\quad\text{for some }\sigma\in C^\infty(M),\\
&F^i_{,r}=0\quad \text{for } 1\leq i\leq n,
\end{split}
\end{equation}
If we choose $g\in[g]$ and $h\in[h]$ suitably such that $f$ is homothetic, that is $\sigma$ constant, then (\ref{eq:1asympisomcondi})-(\ref{eq:3asympisomcondi}) further imply at $r=0$
\begin{align}
\label{eq:higherderivativeconditionforFiandFs}
\begin{split}
     &F^{ i}_{,rr}=F^{ i}_{,rrr}=0\quad\text{for } 1\leq i\leq n,\\
    &F^s_{,rr}=0,
\end{split}
\end{align}
and at $r=0$
\begin{align}\label{eq:morehigherderivativeconditionforFs}
    F^s_{,rrr}=0.
\end{align}
Conversely, if $F$ satisfies (\ref{eq:1stderivativeconditionforF})-(\ref{eq:morehigherderivativeconditionforFs}), then $F$ is an asymptotic local isometry. 
\\

Let $\Sigma\subset M_+$ be a surface orthogonal $M$ and $(t,\lambda)$ its asymptotic isothermal coordinate (\ref{eq:expofrandx}). Since we are considering (\ref{eq:1stderivativeconditionforF}) to (\ref{eq:morehigherderivativeconditionforFs}), it's better to introduce a change of variables on $\Sigma$, $(t,\lambda)\mapsto (t,\ubar r(t,\lambda))$, where $\ubar r(t,\lambda)$ is equal to the right-hand side of $r(t,\lambda)$ in (\ref{eq:expofrandx}) modulo $O(\lambda^4)$. Then, the expansions of $x^i$ and $r$ in terms of $\ubar r$ are
    \begin{align}
    \label{eq:expofrandxinspecialcoord}
        \begin{matrix*}[l]
        x^i(t,\ubar r)&=&\gamma^i(t)+0+\frac{v^i}{2}\ubar r^2+\frac{u^i}{3|\dot\gamma|^3}\ubar r^3+O(\ubar r^4),
        \\
        r(t,\ubar r)&=&\ubar r+O(\ubar r^4)
        \end{matrix*}
    \end{align}
where $v^i$ and $u^i$ remain the same conditions as in Proposition \ref{prop:expofrandxinisothermal} and in Proposition \ref{prop:expof2ndfundform}. We call $(t,\ubar r)$ the \textit{adapted coordinate} of $\Sigma$.

    \begin{proposition}\label{prop:preserveorthogsurf}
    Let $F\colon M_+\to N_+$ be a local diffeomorphism in terms of (\ref{eq:coorddescripofF}). Assume $F$ maps non-null vectors in $TM$ to non-null vectors in $TN$. Then, $\partial_r\mathscr F=0$ at $r=0$ if and only if $F(\Sigma)$ is orthogonal to $N$ for any surface $\Sigma$ in $\overline{M_+}$ that is orthogonal to $M$ and intersects $M$ along a non-null curve. In addition, if $F$ maps proper surfaces in $\overline{M_+}$ to proper surfaces in $\overline{N_+}$, then $\partial_r\mathscr F=0$ at $r=0$.
    \end{proposition}
    \begin{proof}
    Let $\gamma(t)$ be a non-null curve in a coordinate open set $\mathcal W\subset M$ with the coordinate $\gamma^i(t)$. Then, the adapted coordinate $(t,\ubar r )$ in (\ref{eq:expofrandxinspecialcoord}) locally defines a surface $\Sigma \subset \overline{M_+}$ orthogonal to $M$ with $\Sigma\cap M=\gamma$. Since $F(\Sigma)$ is orthogonal to $N$, then
    \begin{align}\label{eq:vectororthonal}
        dF(X_{\ubar r})-\frac{\langle dF(X_{\ubar r}),dF(X_t)\rangle_{\bar h}}{\langle dF(X_t),dF(X_t)\rangle_{\bar h}} dF(X_t)
    \end{align}
    is orthogonal to $TN\subset TM_+$ where $\{X_t,X_{\ubar r}\}$ is the coordinate basis for the adapted coordinate $(t,\ubar r)$. The orthogonal condition gives at $r=0$
    \begin{align}\label{eq:vectororthogonal2}
        0=F^k_{,r}-\frac{F^i_{,r}\dot\gamma^jh_{ij}}{\langle\dot\gamma,\dot\gamma\rangle_h}\dot\gamma^k.
    \end{align}
    Since $\dot\gamma$ is an arbitrary non-null vector at $t=0$, we get $F^k_{,r}=0$ at $r=0$. Conversely, let $(t,\ubar r)$ be the adapted coordinate of $\Sigma$. Projecting (\ref{eq:vectororthonal}) to $TN\subset T\overline{N_+}$ orthogonally, it gives the right-hand side of (\ref{eq:vectororthogonal2}) which turns to be $0$ due to $F^i_{,r}=0$ at $r=0$.

    If $\gamma(t)$ is a non-null conformal circle, then the formula of the adapted coordinate locally extends $\gamma(t)$ to a proper surface in $\overline{M_+}$. One can follow the same arguments just made to get $\partial_r\mathscr F=0$ at $r=0$.
    \end{proof}
As mentioned at the end of \textsection \ref{sec:confcirc}, the idea for proving Theorem \ref{thm:holocircle2} is to consider a suitable family of proper surfaces $\Sigma_\epsilon$. Then, the dependence of $\epsilon$ in the second fundamental forms of proper surfaces $F(\Sigma_\epsilon)$ may imply $F$ is an asymptotic local isometry where its local conditions are (\ref{eq:1stderivativeconditionforF})-(\ref{eq:morehigherderivativeconditionforFs}). However, recalling Proposition \ref{prop:expofrandxinisothermal}, Proposition \ref{prop:expof2ndfundform} and the adapted coordinate (\ref{eq:expofrandxinspecialcoord}), a proper surface $\Sigma$ is characterized by $v^i$, $u^i$ and $\gamma$ being an unparametrized conformal circle. Therefore, we can utilize the adapted coordinate of $F(\Sigma_\epsilon)$ to avoid the tedious computation of the second fundamental forms. The following Lemma \ref{lem:coordchange} and Proposition \ref{prop:holocircle} respectively give the coordinate change of $F(\Sigma)$ to its adapted coordinate and provide that $F$ satisfies (\ref{eq:1stderivativeconditionforF}) and (\ref{eq:higherderivativeconditionforFiandFs}), except for (\ref{eq:morehigherderivativeconditionforFs}), for $F$ preserving proper surfaces.
    \begin{lemma}\label{lem:coordchange}
    Let $(M_+,g_+)$ be a Poincaré-Einstein space in the normal form relative to $g\in[g]$ and $\tilde\Sigma\subset M_+$ be a surface orthogonal to $M$ with $\tilde\Sigma\cap M=\gamma$ being a non-null curve. Assume it has a parametrization $(t,\ubar r)\mapsto (\tilde x^i(t,\ubar r),\tilde r(t,\ubar r))$ from $I\times I$ to $\tilde\Sigma\subset M_+$ where $\tilde x^i_{,\ubar r}=0$ and $\tilde r=0$ both at $\ubar r=0$. The existence of the adapted coordinate $(\eta,\rho)$ of $\tilde\Sigma$ implies there is a coordinate change $t=t(\eta,\rho)$, $\ubar r=\ubar r(\eta,\rho)$ with $t(\eta,0)=\eta$. The coordinate change is in the following modular higher orders.
    \begin{align}
        \begin{matrix*}[l]\label{eq:coordchangeinlemma}
            t(\eta,\rho)&=&\eta+t_{(2)}\rho^2+t_{(3)}\rho^3,\\
            \ubar r(\eta,\rho)&=&r_{(1)}\rho+r_{(2)}\rho^2+r_{(3)}\rho^3
        \end{matrix*}
    \end{align}
    where 
    \begin{align*}\small
        r_{(1)}&=\frac{1}{\tilde r_{,\ubar r}},\;r_{(2)}=-\frac{\tilde r_{,\ubar r\ubar r}}{2(\tilde r_{,\ubar r})^3},\; r_{(3)}=-\frac{1}{6r_{,\ubar r}}\left(2\tilde r_{,t\ubar r}\,r_{(1)}t_{(2)}+\tilde r_{,\ubar r\ubar r\ubar r}(r_{(1)})^3+6\tilde r_{,\ubar r\ubar r} r_{(1)}r_{(2)}\right),
        \\
        t_{(2)}&=\frac{1}{2\langle\dot{\gamma},\dot{\gamma}\rangle}\left(v^i\dot{\gamma}_i- \tilde x^i_{,\ubar r\ubar r}\dot{\gamma}_i(r_{(1)})^2\right), \;t_{(3)}=\frac{1}{6\langle\dot{\gamma},\dot{\gamma}\rangle}\left((r_{(1)}^3) \tilde x^i_{,\ubar r\ubar r\ubar r}\dot{\gamma}_i+6r_{(1)}r_{(2)}\tilde x^i_{,\ubar r\ubar r}\dot{\gamma}_i\right).
    \end{align*}
    The partial derivatives of $\tilde x^i$ and $\tilde r$ above are at $\ubar r=0$. The term $v$ in $t_{(2)}$ is $\nabla_{\dot{\gamma}}\left(\frac{\dot{\gamma}}{\langle \dot{\gamma},\dot{\gamma}\rangle}\right)$. Note that $\dot\gamma^i=\tilde x^i_{,t}(t,0)=\tilde x^i_{,\eta}(\eta,0)$ and $\tilde x^i_{,tt}(t,0)=\tilde x^i_{,\eta\eta}(\eta,0)$.
    \end{lemma}
    
    \begin{proposition}\label{prop:holocircle}
    Let $F\colon M_+\to N_+$ be a local diffeomorphism such that 
    \begin{enumerate}[label=(\roman*)]
    \item $dF(v)$ is non-null (resp. null) $\;\forall$non-null (resp. null) $v\in TM$ if $p\neq q$;
    \item sgn$\langle v,v\rangle_g$=\;sgn$\langle dF(v),dF(v)\rangle_h$ $\forall v\in TM$ if $p=q$.
    \end{enumerate}
    The $F$ satisfies (\ref{eq:1stderivativeconditionforF}) and (\ref{eq:higherderivativeconditionforFiandFs}) if and only if $F$ maps proper surfaces in $\overline{M_+}$ to proper surfaces in $\overline{N_+}$.
    \end{proposition}
\begin{proof}
    Let $\gamma(t)$ be a non-null conformal circle in a coordinate open set of $M$ with the initial conditions $\dot\gamma^k_0$ and $\ddot\gamma^k_0$ at $\gamma(0)=p$. Let $\Sigma\subset \overline{M_+}$ be its extended proper surface defined by the adapted coordinate $(t,\ubar r)$ in (\ref{eq:expofrandxinspecialcoord}). Since $F(\Sigma)\cap N$ is still an unparametrized conformal circle from Proposition \ref{prop:expof2ndfundform}, we know $f\colon M\to N$ is a conformal local diffeomorphism due to Theorem \ref{thm:conformalinftycircle}. Without loss of generality, we assume that $f$ is the identity map on $(M,g)$ where $g\in[g]$ and that the local coordinate of its extended local diffeomorphism is $F\colon (x^i,r)\mapsto (F^i(x,r),F^r(x,r))$ where $r$ is the geodesic defining function of $g$. From Proposition \ref{prop:preserveorthogsurf} and Lemma \ref{lem:coordchange}, we let $(\eta(t,\ubar r),\rho(t,\ubar r))$ be the coordinate change of $F(\Sigma)$ to its adapted coordinate.
    
    Let $\dot\gamma^k_0=V^k$, $\ddot\gamma^k_0=\epsilon A^k$ with $|V| =1$ and $\langle V,A\rangle=0$ where $\epsilon\in\mathbb R$ is arbitrary near $1$. The variable $\epsilon$ gives an $\epsilon$-family of non-null conformal circles $\gamma_\epsilon(t)$. The extended proper surfaces $\Sigma_\epsilon$ is defined by the formula of the adapted coordinate (\ref{eq:expofrandxinspecialcoord}) where the coefficients in the expansion of $x^i(t,\ubar r)$ depend on $\gamma_\epsilon(t)$. Because $\Sigma_\epsilon$ depends smoothly on $\epsilon$, we know the adapted coordinate for $F(\Sigma_\epsilon)$ defined from Lemma \ref{lem:coordchange} depends smoothly on $\epsilon$. Since $F(\Sigma_\epsilon)$ are proper surfaces, we have 
    \begin{align}\label{eq:coordchangev}
 \nabla_{\dot\gamma_\epsilon}\left(\frac{\dot\gamma_\epsilon}{\langle\dot\gamma_\epsilon,\dot\gamma_\epsilon\rangle}\right)^i=F^i_{,\rho\rho} \qquad \text{at } \rho=0
    \end{align}
    and
    \begin{align}\label{eq:coordchangen}
        0=F^i_{,\rho\rho\rho} \qquad \text{at }\rho=0.
    \end{align}
Following the formulas and the conventions from (\ref{eq:coordchangeinlemma}), we know $t_{(2)}=-\frac{V_jF^j_{,rr}}{(F^r_{,r})^2}$ at the point $p$. Using chain rule on $F^i_{,\rho\rho}$ for $(t(\eta,\rho),\ubar r(\eta,\rho))$, the straightforward computation for (\ref{eq:coordchangev}) at $p$ is
\begin{align}
\label{eq:pf_chainruleonv}
    \epsilon A^i=\frac{1}{(F^r_{,r})^2}\left(\epsilon A^i+F^i_{,rr}-V_jF^j_{,rr}\,V^i\right). 
\end{align}
Therefore, $F^r_{,r}=1$ and $F^i_{,rr}=0$ at $r=0$. The results we got imply $r_{(2)}=-F^r_{,rr}$ at $r=0$ and $t_{(3)}=\frac{1}{6}F^j_{,rrr}V_j$ at $p$ from (\ref{eq:coordchangeinlemma}). Doing chain rule again, (\ref{eq:coordchangen}) is equal to the following at $p$.
\begin{align}
    0=-3\epsilon\, F^r_{,rr}\,A^i+F^i_{,rrr}+F^j_{,rrr}V_j \,V^i.
\end{align}
So, $F^r_{,rr}=0$ and $F^i_{,rrr}=0$ at $r=0$.

Conversely, assume $(F^i(x,r),F^r(x,r))=(x^i+O(r^4),r+O(r^3))$ where we assume $f$ is the identity map on $(M,g)$. Let $\Sigma\subset \overline{M_+}$ be a proper surface and $(t,\ubar r)$ be its adapted coordinate. From Proposition \ref{prop:preserveorthogsurf}, we know $F(\Sigma)$ orthogonal to $N$. Considering the coordinate change $(\eta(t,\ubar r),\rho(t,\ubar r))$ of $F(\Sigma)$ to its adapted coordinate and following the formula from (\ref{eq:coordchangeinlemma}), we have
\begin{align*}
    r_{(1)}=1,\;r_{(2)}=0,\;t_{(2)}=0,\;t_{(3)}=0.
\end{align*}
Computing $F^i_{,\rho\rho}$ and $F^i_{,\rho\rho\rho}$ directly from chain rule, we have at $\rho=0$
\begin{align*}
    F^i_{,\rho\rho}=\nabla_{\dot\gamma}\left(\frac{\dot\gamma}{\langle\dot\gamma,\dot\gamma\rangle}\right)^i,\; F^i_{,\rho\rho\rho}=0
\end{align*}
which implies $F(\Sigma)$ is proper by (\ref{eq:expofrandxinspecialcoord}) and Proposition \ref{prop:expof2ndfundform}.
\end{proof}

    \begin{theorem}
    \label{thm:holocircle2}
    Let $(M_+,g_+)$ and $(N_+,h_+)$ be Poincaré-Einstein manifolds for $(M,[g])$ and $(N,[h])$ respectively with same signature $(p+1,q)$. Given a local diffeomorphism $F\colon M_+\to N_+$ such that it smoothly extends a local diffeomorphism $f\colon M\to N$. Assume $F$ satisfies
    \begin{enumerate}[label=(\roman*)]
    \item $dF(v)$ is non-null (resp. null) $\;\forall$non-null (resp. null) $v\in TM$ if $p\neq q$,
    \item sgn$\langle v,v\rangle_g$=\;sgn$\langle dF(v),dF(v)\rangle_h$ $\forall v\in TM$ if $p=q$.
    \end{enumerate}
    If the $F$ maps proper surfaces in $M_+$ to proper surfaces in $N_+$, then there is a local diffeomorphism $G$ on an open neighborhood $W$ of $M\subset \overline{M_+}$, $$G\colon W\to M_+,$$
    where $G$ smoothly extends the identity map on $M$, such that $\tilde F=F\circ G$ is an asymptotic local isometry.
    \end{theorem}

\begin{proof}
    From Proposition \ref{prop:holocircle}, we can choose $g\in[g]$ and $h\in[h]$ suitably to let $f$ be a local isometry. Consider the identification of $F$ in (\ref{eq:coorddescripofF}) 
    \begin{align*}
        F\colon \mathcal U\to \mathcal V,\quad (x,r)\mapsto(\mathscr F(x,r),F^s(x,r)).
    \end{align*}
    We aim to find out an open set $\mathcal W\subset \mathcal U$ such that it contains $M\times\{0\}$ and the following map is well-defined
    \begin{align}
    \label{eq:thm_defofG}
        \begin{matrix*}[l]
            G\colon &\quad\mathcal W &\to& \quad\mathcal U\\
        &(x,r)&\mapsto& (x,r-\mathscr R(x)r^3)=(x,p_x(r))
        \end{matrix*}
    \end{align}
    where $\mathscr R(x)=\partial^3_r F^s(x,0)/6$ and $p_x(r)=r-\mathscr R(x)r^3$. 
    
    For any $x\in M$, it has open neighborhoods $\mathcal B_x,\,\mathcal N_x$ in $M$ such that
    $$
    x\in \mathcal B_x\subset\subset\mathcal N_x\quad\text{and }\mathcal N_x\times[0,\epsilon)\subset \mathcal U\quad\exists\, \epsilon>0
    $$
    where $\mathcal B_x\subset\subset\mathcal N_x$ means its closure $\overline{\mathcal B_x}$ is compact in $\mathcal N_x$. Here we choose $\epsilon$ small enough such that the polynomial $p_y(r)$ of $r$ on $[0,\epsilon)$ is strictly increasing for all $y\in\overline{\mathcal B_x}$. Hence, there exists $0<\epsilon'_x<\epsilon$ so that 
    \begin{align*}
    \mathcal B_x\times [0,\epsilon'_x) &\to \mathcal B_x\times[0,\epsilon)\subset \mathcal U    \\
    (y,r)&\mapsto (y,p_y(r)).
    \end{align*}
    Let $\mathcal W= \bigcup_{x\in M}\mathcal B_x\times[0,\epsilon'_x)$. Since we know the asymptotic expansion of $F$ from Proposition \ref{prop:holocircle}, we get $\tilde F=F\circ G \colon\mathcal W\to \mathcal V$ is an asymptotic local isometry.  
\end{proof}

    \begin{corollary}
    \label{cor:partial3 F}
    Let $F\colon M_+\to N_+$ and $G\colon W\to M_+$ be local diffeomorphisms as stated in Theorem \ref{thm:holocircle2}. Assume there is a geodesic defining function $r$ for some $g\in [g]$ and $C\geq 0$ such that
    \begin{equation}
        |\partial^3_r(s\circ F)(p)|\leq C\quad \forall p\in M\subset M_+,
    \end{equation}
    where $\partial_r={}^{\bar g}\nabla r$ is the gradient of $r$ with respect to $\bar g=r^2g_+$ and $s$ is the geodesic defining function of $h\in[h]$ to make $f\colon M\to N$ be a local isometry.
    \begin{enumerate}[label=(\roman*)]
        \item If $C>0$, then $G$ can be chosen as an embedding with its image while $W$ is small enough. \label{item1:cor:partial3 F}
        \item If $C=0$, then $G$ can be chosen to be the identity map on $W$. Particularly, $F$ is an asymptotic local isometry.\label{item2:cor:partial3 F}
    \end{enumerate}
    \end{corollary}
    \begin{proof}
        Recall the definition of $G$ in (\ref{eq:thm_defofG}),
        \begin{align*}
            \begin{matrix*}[l]
            G\colon &\quad\mathcal W &\to& \quad\mathcal U\\
            &(x,r)&\mapsto& (x,r-\mathscr R(x)r^3)=(x,p_x(r)).
        \end{matrix*}
        \end{align*}
        Let $\mathcal W$ be small enough so that $\mathcal W\subset M\times [0,\epsilon)$ for some $\epsilon>0$. If $\epsilon$ is small enough, then $p_x(r)$ is strictly increasing because for $r\in [0,\epsilon)$
        \begin{align*}
            \partial_r p_x(r)=1-3\mathscr R(x)r^2>1-3Cr^2.
        \end{align*}
        Hence, $G$ is an open injective immersion for $C>0$. The case for $C=0$ is straightforward due to the definition of $G$.
    \end{proof}
    
    The following proposition gives geometric conditions to satisfy the presumptions of Corollary \ref{cor:partial3 F}.
    \begin{proposition}
    \label{prop:(h_+-g_+)_rr}
    Let $F\colon M_+\to N_+$ and $G\colon W\to M_+$ be local diffeomorphisms as stated in Theorem \ref{thm:holocircle2}. Assume there is a geodesic defining function $r$ for some $g\in [g]$ and $a\geq 0$ such that
    \begin{equation}
        \label{eq:prop_geomcondiofF}
        (F^*h_+-g_+)({}^{\bar g}\nabla r,{}^{\bar g}\nabla r)=O(r^a).
    \end{equation}
    Let $s$ be a geodesic defining function of $h\in[h]$ to make $f$ be a local isometry.
    \begin{enumerate}[label=(\roman*)]
        \item If $a=0$, there exists $C>0$ such that $|\partial^3_r(s\circ F)|\leq C$ on $M$.
        \label{item1:prop:(h_+-g_+)_rr}
        \item If $a=1$, we have $\partial^3_r(s\circ F)=0$ on $M$.
        \label{item2:prop:(h_+-g_+)_rr}
    \end{enumerate}
    \end{proposition}
    \begin{proof}
    Consider the identification of $F$ near $M\subset\overline{M_+}$ and $N\subset\overline{N_+}$ in (\ref{eq:coorddescripofF}), $\mathcal U\to \mathcal V,\;(x,r)\mapsto(\mathscr F(x,r),F^s(x,r))$. Then, (\ref{eq:prop_geomcondiofF}) is equivalent to 
    \begin{align*}
        O(r^{a+2})=F^s_{,r}F^s_{,r}+F^i_{,r}F^j_{,r}(h_s\circ F)_{ij}-(F^s/r)^2,
    \end{align*}
    where the right-hand side above is exactly from (\ref{eq:2asympisomcondi}) while considering local conditions of an asymptotic local isometry. Then, the Taylor expansion for the right-hand side gives the following.
\begin{align*} 
    O(r^{a+2})=&\sum^2_{c=0}\left\{\sum_{c\geq b\geq 0} \frac{(b+1)(c-b+1)}{(b+1)!\,(c-b+1)!}F^s_{(b+1)}F^s_{(c-b+1)}\right.
    \\
    &+\sum_{c\geq b+d\geq 0}\frac{(b+1)(c-b-d+1)}{(b+1)!\,(c-b-d+1)!\,d!}F^i_{(b+1)}F^j_{(c-b-d+1)}(h_{(d)})_{ij}
    \\
    &\left.-\sum_{c\geq b\geq 0}\frac{1}{(b+1)!\,(c-b+1)!}F^s_{(b+1)}F^s_{(c-b+1)}\right\}r^c,
\end{align*}
where $F^s_{(b)}$, $F^i_{(b)}$ and $h_{(b)}$ mean the $b$th-order partial derivative of $r$ at $r=0$ of $F^s$, $F^i$ and $h_s\circ F$ respectively. So, we have from above
\begin{equation}
\label{eq:condiof(h_+-g_+)(r,r)fromtaylor}
    \begin{split}
    0&= F^i_{,r}F^j_{,r}h_{ij}\;,
    \\
    0&= F^s_{,r}F^s_{,rr}+3F^i_{,r}F^j_{,rr}h_{ij}\;,
    \\
    O(r^a)&= \frac{2}{3}F^s_{,r}F^s_{,rrr}+\frac{3}{4}(F^s_{,rr})^2+F^i_{,rr}F^j_{,rr}h_{ij},  
    \end{split}
\end{equation}
where the third equality above is when $F^i_{,r}=0$ at $r=0$. We know $F^s_{,r}=1$, $F^s_{,rr}=0$ and $F^i_{,rr}=0$ at $r=0$ from Proposition \ref{prop:holocircle}. This completes the proof.
\end{proof}

\begin{remark}
        It is straightforward to observe $F^i_{,r}=0$ and $F^s_{,rr}=0$ at $r=0$ from (\ref{eq:condiof(h_+-g_+)(r,r)fromtaylor}) when considering Riemannian conformal classes $[g]$ and $[h]$. However, (\ref{eq:condiof(h_+-g_+)(r,r)fromtaylor}) alone still can't simply imply $F^s_{,rrr}= 0$ at $r=0$.
\end{remark}

\bibliographystyle{plain}
\bibliography{bibliography}
\end{document}